\documentclass[12pt]{article}
\newtheorem{thm}{Theorem}
\newtheorem{cor}[thm]{Corollary}
\newtheorem{example}[thm]{Example}
\newtheorem{definition}[thm]{Definition}
\begin{document}
%\renewcommand{\baselinestretch}{1.6}

%\pagestyle{plain}
%\tightenlines
%\begin{document}
%\renewcommand{\baselinestretch}{1.6}
%\preprint{}
%\begin{frontmatter}
\title{Less Regular Exceptional and Repeating Prime Number Multiplets}
\author{H. J. Weber\\Department of Physics\\
University of Virginia\\Charlottesville, 
VA 22904, U.S.A.}
%\date{\today}
\maketitle
\begin{abstract}
New exceptional (i.e. non-repeating) prime number 
multiplets are given and formulated in terms of 
arithmetic progressions, along with laws governing 
them. Accompanying repeating prime number multiplets 
are pointed out. Prime number multiplets with less 
regular distances are studied.    
\end{abstract}
%\vspace{3ex}
%\pacs{}
\leftline{MSC: 11N05, 11N32, 11N80}
\leftline{Keywords: Exceptional, repetitious 
prime number multiplets.} 

%\vspace{3ex}

\section{Introduction}

In Refs.~\cite{hjw1},\cite{hjw2},\cite{hjw3} 
a variety of prime number multiplets have been 
discussed, most of which exhibit a regular 
distance pattern. The reason for restricting 
attention to them is the enormously complex  
mix of regularities with chaotic properties 
of prime numbers. 

The well-known exceptional triplet $3,5,7$ 
is the only case of three successive primes 
in the arithmetic progression $3+2n$ for 
$n=0,1,2.$ All others are composites or single 
primes, such as $23$ for $n=10,$ and $37$ 
for $n=17,$ and ordinary twins like $11,13$ 
for $n=4,5$ and $17,19$ for $n=7,8$ etc.    

More general exceptional triplets~\cite{hjw1},
\cite{hjw2} such as $3,3+2d,3+4d$ of primes at 
equal distance $2d,~(3,d)=1$ translate into 
successive prime values in the arithmetic 
progression $3+2dn.$ There is at most one prime 
number triplet in it (for $n=0,1,2$) and all 
others are composites, single primes or twins 
at the distance $2d,$ but no $k-$tuple  of 
primes for $k>3.$        

For any odd prime $p,$ in terms of the 
arithmetic progression $p+2dn$ with $(d,p)=1,$ 
there is at most one $p-$tuple of primes for 
$n=0,1,2,\ldots,p-1$~\cite{hjw2}, but no 
$k-$tuples for $k>p.$ In Section~3 it is 
shown that there are accompanying 
$(p-1)-$tuples of primes that usually repeat.  

We also generalize the prime $3$ and triplets 
in (ii) of Theor.~2.3~\cite{hjw1} to an 
arbitrary odd prime $p$ and $p-$tuples of 
primes. This continues the quest for 
uncovering more and deeper laws governing 
prime number multiplets.  

\section{New Exceptional Prime Number\\Multiplets}

Let us start with a few samples of prime number 
multiplets that generalize the exceptional 
triplets in (ii) of Theor.~2.3~\cite{hjw1}. 
  
Given any distance $2d_1$ not divisible by $5$ 
such that $5+2d_1$ is prime, we can extend  
$5, 5+2d_1,5+2d_2,5+2d_3,5+2d_4$ to a 
maximum-length exceptional quintet of prime 
numbers using $5|d_j-d_1$. That is to say, 
there is at most one such quintet for a given 
distance pattern, and they cannot repeat.  
 
{\bf Example~2.1.} For $d_1=1,~3,~4,~6,~7,\dots$ 
\begin{eqnarray*}
&&5,~7,~19,~31~43~[2,2(1+5),2,2(1+5)]-7,-19,-31,-43\\
&&5,~11,~17,~23,~29~[6,6,6,6]\\
&&5,~11,~37,~43,~59~[6,2(3+10),6,2(3+5)]-11,-17,-23,-29\\
&&5,~13,~31,~59,~67~[8,2(4+5),2(4+10),8]-13,-31,-59,-67\\
&&5,~17,~29,~41,~53~[12,12,12,12]-7,-19,-31,-43.
\end{eqnarray*}
The continuation to negative numbers is given behind  
the distance pattern in brackets. 

None of these quintets (nonets) can be extended to $6$ 
(or $10$) primes in a row with the same or similar 
distance pattern. And they do not repeat, i.e. are all 
exceptional. The same is the case for $5,19,23,37,41$ 
at distances $[14,4,14,4]$ with $5|14-4,~5\not| 14.$ 

Similar rules apply to septets for $d_1=2,3,5,6,\ldots$.
\begin{eqnarray*}
&&7,11,29,47,79,83,101~[4,2(2+7),18,2(2+2\cdot7),4,18]
\\&&3,-29,-47,-79,-83,-101\\
&&7,13,19,53,59,79,113~[6,6,2(3+2\cdot 7),6,2(3+7),34]\\
&&7,17,41,79,89,113,137~[10,2(5+7),2(5+2\cdot 7),10,24,24]
\\&&7,19,31,43,83,137,149~[12,12,12,2(6+2\cdot 7),2(6+
3\cdot 7),12].
\end{eqnarray*}
These are all special cases of the following general 
law.     

{\bf Theorem~2.2.} {\it Let $p$ be an odd prime, 
$p|d_j-d_1, p\not| d_1, d_j>0.$ Then there is at most 
one prime number $p-$tuple $p, p+2d_1,\ldots,p+2d_{p-1}$ 
with distance pattern $[2d_1,2d_2,\ldots,2d_{p-1}]$.} 

{\it Proof.} If $p, p+2d_1,\ldots,p+2d_{p-1}$ is 
a prime number p-tuple, then $p_3=p_2+2d_2\equiv 
p_1+4d_1\pmod{p},\ldots,p_{p-1}\equiv p_1+2(p-1)
d_1\pmod{p}.$ There is no other such $p-$tuple,  
because one of $p$ odd numbers in a row at the  
same distance $2d_1$ (mod $p$) is divisible by 
$p.~\diamond$ 

{\bf Corollary~2.3.} {\it There is a $p-$tuple 
of primes with a distance pattern $2d_j\pmod{p}$  
of Theor.~2.2.}

{\it Proof.} With $pn+2d_1$ forming an arithmetic 
progression, Dirichlet's theorems allows us 
picking $n_1$ so that $pn_1+2d_1=p_1$ is prime. 
Likewise, we pick $n_2$ in $pn_2+2d_2=p_2$ so 
it is prime, etc. The proof in Theor.~2.2 shows 
that there is no other such $p-$tuple of primes.
~$\diamond$

Using this principle we can construct exceptional 
$p-$tuples of primes for any odd prime number $p$ 
as follows. 

{\bf Example~2.4.} The exceptional septet $7,13,
47,67,73,79,113$ has the distances $[6,2(3+2\cdot 7),
2(3+7),6,6,34],$ where each distance has the form 
$2\cdot 3+2\cdot 7n$ in accord with Theor.~2.2. It 
cannot be continued to an octet because 
$113+2(7n+3)=7(17+2n)$ factorizes. 

Likewise, the $11-$tuple of primes  
\begin{eqnarray}
11,13,37,61,107,109,199,223,269,271,383
\label{11tuple}
\end{eqnarray}
with distances 
\begin{eqnarray}
[2,2(1+11),24,2(1+2\cdot 11),2,2(1+4\cdot 11),24,
46,2,2(1+5\cdot 11)] 
\label{1}
\end{eqnarray}
cannot be continued to a $12-$tuple, as 
$383+2(11n+1)=\\11(5\cdot 7+2n)$ factorizes. 

The exceptional $13-$tuple of primes 
\begin{eqnarray}
&&13,17,47,103,107,137,167,197,227,257,313,317,347 
\end{eqnarray}
with the distance pattern 
\begin{eqnarray}
[4,2(2+13),2(2+2\cdot 13),4,30,30,30,30,30,56,4,30] 
\end{eqnarray}
cannot be continued because 
$347+2(2+13n)=13(27+2n)$ factorizes. 

The $17-$tuple of primes 
\begin{eqnarray}\nonumber
&&17,~19,~89,~193,~229,~367,~607,~643,~883,
~919,~1193,~1229,\\&&1231,~1301,~1303,~1373,
~1409
\end{eqnarray}
with the distances 
\begin{eqnarray}\nonumber
&&[2,2(1+2\cdot 17),2(1+3\cdot 17),2(1+17),
2(1+4\cdot 17),\\&&2(1+7\cdot 17),36,240,36,
2(1+8\cdot 17),36,2,70,2,70,36]
\end{eqnarray}
stops because $1409+2(1+17n)=17(83+2n)$ 
factorizes. These cases follow the general 
factorization 
\begin{eqnarray}
p+2d_1(p-1)+(2d_1+np)=p
(2d_1+n).
\end{eqnarray}

In a tour de force, we give the following 
$43-$tuple that goes two steps beyond Euler's 
optimal sequence of $41$ primes generated by 
$x(x-1)+41,~x=0,1,\ldots,40:$
\begin{eqnarray*}
&&43,47,137,227,317,751,1013,1103,1193,1283,1373,1549,\\
&&1553,1901,2593,2683,2687,2777,2953,2957,3391,3739,\\
&&4001,4091,4783,4787,4877,4967,5573,5749,5839,\\
&&5843,6277,6367,7489,8009,8443,8447,8537,8627,8803,\\
&&8807,9241
\end{eqnarray*} 
with the distance pattern
\begin{eqnarray*}
&&[4,2(2+43),2(2+43),2(2+43),2(2+5\cdot 43),
2(2+3\cdot 43),\\&&2(2+43),2(2+43),2(2+43),
2(2+43),2(2+2\cdot 43),4,\\&&2(2+4\cdot 43),
2(2+8\cdot 43),2(2+43),4,2(2+43),2(2+2\cdot 43),
4,\\&&2(2+5\cdot 43),2(2+4\cdot 43),2(2+3\cdot 43),
2(2+43),2(2+8\cdot 43),\\&&4,2(2+43),2(2+43),
2(2+7\cdot 43),2(2+2\cdot 43),2(2+43),4,\\&&
2(2+5\cdot 43),2(2+43),2(2+13\cdot 43),2(2+6
\cdot 43),2(2+5\cdot 43),\\&&4,2(2+43),2(2+43),
2(2+2\cdot 43),4,2(2+5\cdot 43)].
\end{eqnarray*}
In each step, as outlined in the proof of 
Theor.~2.5 below, we pick the next possible 
prime. Yet, sometimes there are gaps of hundreds.  
In general, it is much easier to search for the 
first long prime sequence than uncover record 
setting $p-$tuples at equal distances~\cite{hjw3} 
of comparable length. 

{\bf Theorem~2.5} {\it Given any odd prime 
$p,$ there are infinitely many $2d_1>0$ such 
that $p+2d_1$ is prime and $p-2$ multiples 
$pj_i$ of $p$ so that the sequence $p,p+2d_1,
p+4d_1+2pj_1,\ldots,p+2(p-1)d_1+2pj_{p-2}$ 
forms a $p-$tuple of primes. Each $p-$tuple 
of primes is of maximum length and exceptional, 
i.e. will not repeat.}  

{\bf Remark~2.6.} Since the general multiplet 
member has the form $p(1+2j_k)+2kd_1,$ this 
result may be viewed as the existence of 
maximum-length succession of primes in some 
generalized arithmetic progressions, and the 
following proof clarifies what is meant by this.   

{\it Proof.} Given the odd prime $p,$ we 
pick a $d_1$ so that $p+2d_1=p_1$ is prime. 
There are infinitely many such values by 
Dirichlet's theorem, because $p+2d_1$ is 
an arithmetic progression with $d_1$ (and 
$p\not| d_1$) running. Next, we pick $j_1$ 
in $p+2d_1+(2d_1+2j_1p)=p_2$ so that $p_2$ 
is prime. Again by Dirichlet's theorem 
there is an infinitude of such values $j_1,$ 
because $p(1+2j_1)+4d_1$ is an arithmetic 
progression. And just as in step 1, each 
of these choices leads to a complete 
$p-$tuple of primes, and so on. In step 
$p-1,$ we pick $j_{p-2}$ so that 
$p(1+2j_{p-2})+2(p-1)d_1$ is prime. This 
being an arithmetic progression with 
$j_{p-2}$ running, it can be done again by 
Dirichlet's theorem in infinitely many 
ways. This completes the $p-$tuple and 
the proof, because no further step is 
possible in view of the factorization 
\begin{eqnarray}\nonumber
&&p(1+2j_{p-2})+2(p-1)d_1+(2d_1+2j_{p-1}p)\\
&=&p[2d_1+1+2j_{p-2}+2j_{p-1}].~\diamond
\end{eqnarray}   
This construction testifies to the 
unbelievable variety, richness and 
complexity of the sequence of ordinary 
prime numbers.   

We could also have walked backward at any 
step, as is shown in the following examples. 

{\bf Example~2.7.} The quintet 
\begin{eqnarray}
5,7,19,11,13~[2,2(1+5),2(1-5),2]
\label{5tuple}
\end{eqnarray} 
is stopped by $13+2(1+5n)=5(3+2n);$
\begin{eqnarray}
5,7,19,11,3~[2,2(1+5),2(1-5),2(1-5)]
\end{eqnarray}
by $3+2(1+5n)=5(1+2n);$
\begin{eqnarray}
5,7,-11,-19,-17~[2,2(1-2\cdot 5),2(1-5),2] 
\end{eqnarray}
by $-17+2(1+5n)=5(-3+2n);$
\begin{eqnarray}
5,7,-11,-19,-37~[2,2(1-2\cdot 5),2(1-5),2(1-2\cdot 5)],
\end{eqnarray}
by $-37+2(1+5n)=5(-7+2n).$ Only walking straight left  
\begin{eqnarray}
5,-3,-11,-19,-37~[2(1-5),2(1-5),2(1-5),2(1-2\cdot 5)]
\end{eqnarray}
yields an optimal nonet upon extending the quintet in 
Eq.~\ref{5tuple}. 

{\bf Corollary~2.8.} {\it By working to 
the left, some $p-$tuples of primes may 
be extended to optimal $(2p-1)-$tuples.}

{\it Proof.} Essentially the same proof 
as for Theor.~2.5 works to the left, 
generating $p-$tuples involving negative 
prime numbers.~$\diamond$  

{\bf Example~2.9.} Example~2.1 lists several cases 
of nonets and a $13-$tuple. Equation~\ref{11tuple} 
continues as 
\begin{eqnarray}
11,-13,-37,-61,-107,-109,-199,-223,-269,-271,-317
\end{eqnarray} 
with the distance pattern
\begin{eqnarray}\nonumber
&&[-2(1+11),-24,-24,-2(1+2\cdot 11),-2,
-2(1+4\cdot 11),-24,\\&&-46,-2,46]. 
\end{eqnarray}

{\bf Induction principle for primes.} 
{\it Let $A_1,\ldots,A_n$ be a finite set of 
formulas involving a finite number of primes. 
Let $A_1,\ldots,A_n$ be true for the primes 
$p_1,\ldots,p_k.$ If $A_1,\ldots,A_n$ are 
taken to be true for a general set of primes 
$P_1,\ldots,P_k$ and it is shown that there 
are primes $Q_1>P_1,\ldots,Q_k>P_k$ so that 
$A_1,\ldots,A_n$ are true, then $A_1,\ldots,A_n$ 
hold for an infinitude of prime sets.}       
 
The primes $Q_1,\ldots,Q_k$ may be found by 
Euclid's proof of the infinitude of 
primes~\cite{hw},\cite{pr} or, if arithmetic 
progressions are involved, by using 
Dirichlet's theorem. Cor.~2.3, Theors.~2.5, 
3.4 are applications of the prime number 
induction principle. In all of them $A_j$ 
are arithmetic progressions.    

\section{Repeating Prime Patterns}

The previous section dealt with non-repeating, 
or exceptional, $p-$tuples of primes. If they 
were incredibly numerous and diverse, repeating 
patterns are even more so, as we exemplify in 
this section. 

{\bf Example~3.1.} The exceptional quintet 
\begin{eqnarray}
5,11,17,23,29~[6,6,6,6]
\end{eqnarray}
is followed by a string of quartets at equal 
distance $6:$ 
\begin{eqnarray}\nonumber
&&41,47,53,59;~61,67,73,79;~251,257,263,269;
\\&&601,607,613,619;~641,647,653,659;\ldots
\end{eqnarray} 
Each quartet is preceded and followed by a multiple 
of $5,$ such as $35,65;~55,85;~245,275;\ldots.$
Probably there is an infinitude of such quartets  
but, at fixed equal distance $6,$ this may be as 
hard to prove as any twin prime conjecture. In 
terms of the arithmetic progression $5+6n,$ there 
is just one quintet of primes accompanied by a 
string of quartets of primes, but no $k-$tuples 
for $k>5.$    
 
Similarly, the exceptional septet at equal 
distance~\cite{led} $150$
\begin{eqnarray}
7,157,307,457,607,757,907
\end{eqnarray}
is accompanied by (probably infinitely many) 
$6-$tuples of primes at equal distance $150:$
\begin{eqnarray}\nonumber
&&73,223,373,523,673,823;\\\nonumber
&&2467,2617,2767,2917,3067,3217;\\\nonumber
&&4637,4787,4937,5087,5237,5387;\\\nonumber
&&6079,6229,6379,6529,6679,6829;\\\nonumber
&&7717,7867,8017,8167,8317,8467;\\
&&13163,13313,13463,13613,13763,13913;\ldots
\end{eqnarray} 
Again, each $6-$tuple is preceded and followed, 
at the same distance $150$, by a multiple of $7,$ 
such as $823+150=7\cdot 139,73-150=-7\cdot 11;~
2467-150=7\cdot 331,3217+150=7\cdot 481;~13163-
150=7\cdot 1859,13913+150=7\cdot 2009.$ 

Of course, this also holds for the record setting 
$11-$tuples of primes~\cite{hjw3}. The first 
$11-$tuple at equal distance $1536160080$ 
is followed by the $10-$tuples at the same distance 
\begin{eqnarray}\nonumber
&&2009803217,3545963297,5082123377,6618283457,
\\\nonumber
&&815443537,9690603617,11226763697,12762923777,
\\\nonumber
&&14299083857,15835243937;\\\nonumber
&&2622695717,4158855797,5695015877,7231175957,
\\\nonumber
&&8767336037,10303496117,11839656197,13375816277,
\\\nonumber
&&14911976357,16448136437;\\\nonumber
&&2646083851,4182243931,5718404011,7254564091,
\\\nonumber
&&8790724171,10326884251,11863044331,13399204411,
\\\nonumber
&&14935364491,16471524571;\\\nonumber
&&3117107701,4653267781,6189427861,7725587941,
\\\nonumber
&&9261748021,10797908101,12334068181,13870228261,
\\\nonumber
&&15406388341,16942548421;\\\nonumber
&&3178320413,4714480493,6250640573,7786800653,
\\\nonumber
&&9322960733,10859120813,12395280893,13931440973,
\\\nonumber
&&15467601053,17003761133;\\\nonumber
&&3276952243,4813112323,6349272403,7885432483,
\\\nonumber
&&9421592563,10957752643,12493912723,14030072803,
\\&&15566232883,17102392963;\ldots.
\end{eqnarray} 
Again, each decuplet is preceded and followed by a 
multiple of $11,$ such as 
\begin{eqnarray*}
&&2009803217-1536160080=11\cdot 43058467;\\
&&15835243937+1536160080=11\cdot 1579218547. 
\end{eqnarray*}

{\bf Remark~3.2.} Each exceptional 
$p-$tuple of primes at a given distance pattern 
is accompanied by a (possibly empty, or finite, 
but usually infinite) set of $(p-1)-$tuples of 
primes that are preceded and followed by 
multiples of the prime $p.$ This is how 
$(p-1)-$tuples are kept from extending into 
$p-$tuples.

{\bf Example~3.3.} The quintets of Example~2.1 
are accompanied by the following repeating 
quartets with distance patterns shorter by one
\begin{eqnarray*}
&&17,29,31,43~[12,2,12]~47,59,61,73;\\
&&137,149,151,163;~167,179,181,193;\ldots .
\end{eqnarray*}
None can be continued to a quintet using any 
of the distances $2,12,$ because $5|17-2,43+2;
~5|17-12,43+12;~\ldots$ i.e. they are preceded 
and followed by multiples of $5.$ The first 
prime $p_1\equiv 2\pmod{5},$ the second prime 
$p_2\equiv 4\pmod{5},$ the third $p_3\equiv 1
\pmod{5}$ and the 4th $p_4\equiv 3\pmod{5},$ 
which follows from $2, 12\equiv 2\pmod{5}.$ 
Likewise, 
\begin{eqnarray*}
&&41,47,73,79~[6,26,6]~191,197,223,229;\\
&&11,17,43,59~[6,26,16]~41,47,73,89;\\
&&131,137,163,179;~191,197,223,239
\end{eqnarray*} 
with 
\begin{eqnarray*}
&&41,191,11,41,131,191\equiv 1\pmod{5},\\
&&47,197,17,47,137,197\equiv 2\pmod{5},\\
&&73,223,43,73,163,223\equiv 3\pmod{5},\\
&&79,229,59,89,179,239\equiv 4\pmod{5}
\end{eqnarray*}
from $6, 26\equiv 1\pmod{5}.$
Similar rules hold for the first septet in 
Example~2.1 with the accompanying $6-$tuples   
\begin{eqnarray*}
&&431,449,467,499,503,521~[18,18,32,4,18]\\
&&35081,35099,35117,35149,35153,35171;\ldots  
\end{eqnarray*}
though there are tremendous gaps between them. 
Again, they are bracketed by multiples of $7:$ 
\begin{eqnarray*}
&&431-18=7\cdot 59,521+18=7\cdot 77;~431-4=7\cdot 
61,\\&&521+4=7\cdot 75;~431-32=7\cdot 57,521+32=
7\cdot 79;\\&&35081-18=7\cdot 5009,35171+18=7\cdot 
5027;\ldots 
\end{eqnarray*}
For the distance pattern $[4,18,18,32,4]$ the 
$6-$tuples start after a large gap 
\begin{eqnarray*}
&&50047,50051,50069,50087,50119,50123;197887,
\\&&197891,197909,197927,197959,197963;\ldots
\end{eqnarray*} 
For the exceptional $11-$tuple of primes in 
Eq.~(\ref{1}) the first three repeating 
decuplets of primes are the following:  
\begin{eqnarray*}
&&7989996643,7989996667,7989996691,7989996737,
7989996739,\\&&7989996829,7989996853,7989996899,
7989996901,7989996913;\\&&13291266463,13291266487,
13291266511,13291266557,\\&&13291266559,13291266649,
13291266673,13291266719,\\&&13291266721,13291266733;
\\&&14024111323,14024111347,14024111371,14024111417,
\\&&14024111419,14024111509,14024111533,14024111579,
\\&&14024111581,14024111593
\end{eqnarray*} 
with the first prime of each decuplet $\equiv 
2\pmod{p}$ and the following primes at the same 
distances as the exceptional $p-$tuple in 
Eq.~(\ref{1}). None of the decuplets can be extended 
to an $11-$tuple of primes because, upon adding 
any distance $\equiv 2\pmod{p}$ to the last prime 
of any decuplet yields a multiple of $5$ and 
subtracting $2\pmod{p}$ from the first prime of 
each decuplet gives a multiple of $11.$ 
 
{\bf Theorem~3.4.} {\it Let $p$ be an odd prime 
and $p,\ldots,pk_j+2d_1j,j=1,\ldots,p-1$ be an 
exceptional $p-$tuple of primes. Then there are 
infinitely many $(p-1)-$tuples of primes at the 
same distances $2d_1j\pmod{p}, j=1, 2,\ldots,p-1$ 
as for the exceptional $p-$tuple of primes.}  

{\it Proof.} In terms of the arithmetic 
progressions $2d_1j+k_jp,(d_1,p)\\=1$ the 2nd 
prime of the $p-$tuple and the first of each 
$(p-1)-$tuple are in the arithmetic progression 
$2d_1+k_1p,\dots$ the $(j+1)$th of the $p-$tuple 
and $j$th of the $(p-1)-$tuple are in the arithmetic 
progression $2d_1j+k_jp$ for $j=1,2,\dots,p-1.$ 
Dirichlet's theorem for arithmetic progressions 
allows for an infinity of primes in each of these  
arithmetic progressions. We pick the first prime 
in each $(p-1)-$tuple so that they are in 
increasing order, then the second primes similarly, 
etc. Then the distances within each $(p-1)-$tuple 
are $2d_1j\pmod{p}$ which are the same as for 
the leading $p-$tuple. They all start right after 
a multiple of $p$ at the distance $2d_1\pmod{p}$ 
and end with a multiple of $p.~\diamond$    

When the odd prime $p$ in an arithmetic 
progression $p+6dn$ is replaced by 
$p^l,~l\geq 2,$ then $p$ is no longer available 
as the first prime of an exceptional $p-$tuple, 
leaving only $k-$tuples with $k\leq p-1$. 

{\bf Corollary~3.5.} {\it There is no exceptional 
$p-$tuple of primes in the arithmetic progression} 
$p^l+6dn,~l\geq 2, p\not|d,~n=0,1,\dots$. 

\section{More General Prime Sequences}

{\bf Example~4.1.} The arithmetic progressions
$35+6n,~55+6n$ contain at most quartets of primes
\begin{eqnarray*}
41,47,53,59;~251,257,263,269;~641,647,653,659;\ldots
\end{eqnarray*}
\begin{eqnarray*}
61,67,73,79;~601,607,613,619;\ldots
\end{eqnarray*}
where multiples of $5$ usually precede and end them,
such as $59+6=5\cdot 13,~251-6=5\cdot 7^2,~79+6=
5\cdot 17.$ However, the triplet $97,103,109$ in
$55+6n$ is preceded by $91=7\cdot 13,$ where
$5<7<11,~13>11.$ Alternately, the triplet $271,277,
283$ ends before $289=17^2,$ while $361=19^2$
precedes $367,373,379.$ This does not happen for 
quadruplets. There are no exceptional quintets 
because $35=5\cdot 7,55=5\cdot 11$ are composite.

{\bf Theorem~4.2.} {\it Let $p_j|a,~p_1<p_2<
\ldots <p_A$ be all the odd prime divisors of $a$ 
in the arithmetic progression $a+6dn,~(a,6d)=1.$ 
Then multiples of primes $P>p_A$ may eliminate 
$(p_1-1)-$tuples creating $k-$tuple fragments of 
primes starting after a multiple of $p_1$ and 
ending before a multiple of $P,$ or starting 
after a multiple of $P$ and ending before a 
multiple of $p_1,$ where $k\leq p_1-2.$ Prime 
divisors $p_j>p_1$ may play the same role.     

If $q|d,~q>3$ is an odd prime divsor of the distance 
$6d$ then there are no $q-$tuples in $a+6dn.$ Let 
$p'|d$ for all $p'<p_M\leq p_1.$ Then 
$(p_M-1)-$tuples of primes are the longest that 
can occur in $a+6dn$. If $p_M>p_1,$ then $(p_1-1)-$
tuples of primes are the longest in} $a+6dn.$

{\it Proof.} As shown in Cor.~5 of Ref.~\cite{hjw3},
a prime divisor $q|d$ is needed to prevent any
$q-$tuple of primes to occur in $a+6dn.$ Any prime
divisor $p_j|a$ eliminates exceptional $p_j-$tuples
by its multiples and serves to end $k-$tuples or
precedes a $k-$tuple by it. Hence there would be
$(p_j-1)-$tuples of primes if it were not for
multiples of $p_k\neq p_j,~p_k|a$ between two
multiples of $p_j$ that eliminate them. As a
consequence, there are fragments of $(p_j-1)-$
tuples of primes, that is smaller $m-$tuples
that start after a multiple of $p_j$ and end
before a multiple of $p_k.$ Other fragments
start after a multiple of $p_{k_1}$ and end
before a multiple of $p_{k_2},$ or start
after a multiple of $p_k$ and end before a
multiple of $p_j,$ etc. Example~2.1 shows
these cases. Thus, at most $k-$tuples are
allowed with $k\leq p_M-1.$ Thus, the more 
prime divisors $a$ has the fewer $k-$tuples 
of primes will occur in the arithmetic
progression $a+6dn.$ The longer the product of
successive odd prime divisors the distance $2d$
has starting form $3,$ the higher the allowed
$k-$tuples are for $q_M<k\leq p_1-1~\diamond$.

{\bf Corollary~4.3.} {\it Under the conditions 
of Theor.~4.2, there are infinitely many 
$(p_M-1)-$tuples of primes with equal distances 
mod $p_M$ starting after some multiple of $p_1$ 
and ending ahead of some multiple of $p_M.$ If 
$p_M>p_1$ then the multiplets are $(p_1-1)-$
tuples.}  

{\it Proof.} This follows along the lines of 
the proof of Theor.~2.5.~$\diamond$ 

Note carefully that these multiplets of primes 
are not necessarily in the arithmetic progression 
$a+6dn.$

%\end{references}

\begin{thebibliography}{0}  
%\begin{references}

\bibitem{hjw1} Weber, H. J., {\it Regularities 
of twin, triplet and multiplet prime numbers,}  
Global J. of Pure and Applied Math. (2011), 
arXiv:1103.0447 [math.NT].   

\bibitem{hjw2} Weber, H. J., {\it Exceptional 
prime number twins, triplets and multiplets,} 
preprint arXiv:1102.3075 [math.NT]. 

\bibitem{hjw3} Weber, H. J., {\it Remarkable 
and reversible prime number patterns,} preprint  
submitted April 2011.   

\bibitem{hw} Hardy, G. H., Wright, E. M.,   
{\it An Introduction to the Theory of Numbers,} 
Clarendon Press, 5th ed., Oxford (1988). 
  
\bibitem{pr} Ribenboim, P., {\it The New Book 
of Prime Number Records,} Springer, Berlin (1996). 

\bibitem{led} Dickson, L. E., {\it History of 
the Theory of Numbers,} Vol. I, Dover, Mineola, 
N.Y. (2005), p. 426.

\end{thebibliography}
\end{document}